\documentclass[dvips]{proc-l}
\usepackage{epsfig}

\newtheorem{thm}{Theorem}[section]
\newtheorem{cor}[thm]{Corollary}

\theoremstyle{definition}

\theoremstyle{remark}

\numberwithin{equation}{section}

\newcommand{\N}{\mathbb{N}}

\newcommand{\C}{\mathbb{C}}
\newcommand{\ov}{\overline}
\newcommand{\dis}{\displaystyle}
\newcommand{\Notequiv}{/\kern-.6em\hbox{$\equiv$} }

\begin{document}

\title[An inequality for the norm of a polynomial factor]
{An inequality for the norm of a polynomial factor}%
\author{Igor E. Pritsker}%

\address{Department of Mathematics, 401 Mathematical Sciences, Oklahoma State
University, Stillwater, OK 74078-1058, U.S.A.}%
\email{igor@math.okstate.edu}

\thanks{Research supported in part by the National Science Foundation
grants DMS-9970659 and DMS-9707359.}%
\subjclass{Primary 30C10, 30C85; Secondary 11C08, 31A15}%
\keywords{Polynomials, uniform norm, logarithmic capacity, equilibrium
measure, subharmonic function, Fekete points}%

\commby{Albert Baernstein II}%


\begin{abstract}

Let $p(z)$ be a monic polynomial of degree $n$, with complex coefficients, and
let $q(z)$ be its monic factor. We prove an asymptotically sharp inequality of
the form $\|q\|_{E} \le C^n \, \|p\|_E$, where $\|\cdot\|_E$ denotes the sup
norm on a compact set $E$ in the plane. The best constant $C_E$ in this
inequality is found by potential theoretic methods. We also consider
applications of the general result to the cases of a disk and a segment.

\end{abstract}

\maketitle


\section{Introduction}

Let $p(z)$ be a monic polynomial of degree $n$, with complex coefficients.
Suppose that $p(z)$ has a monic factor $q(z)$, so that
$$ p(z)=q(z)\, r(z),$$
where $r(z)$ is also a monic polynomial. Define  the uniform (sup) norm on a
compact set $E$ in the complex plane $\C$ by
\begin{equation} \label{1.1}
\|f\|_{E} := \sup_{z \in E} |f(z)|.
\end{equation}
We study the inequalities of the following form
\begin{equation} \label{1.2}
\|q\|_{E} \le C^n \, \|p\|_E,\quad \deg p = n,
\end{equation}
where the main problem is to find the best (the smallest) constant $C_E$, such
that (\ref{1.2}) is valid for {\it any} monic polynomial $p(z)$ and {\it any}
monic factor $q(z)$.

In the case $E=\ov D$, where $D:=\{z:|z|<1\},$ the inequality (\ref{1.2}) was
considered in a series of papers by Mignotte \cite{Mi}, Granville \cite{Gr}
and Glesser \cite{Gl}, who obtained a number of improvements on the upper bound
for $C_{\ov D}$. D. W. Boyd \cite{Boy1} made the final step here, by proving that
\begin{equation} \label{1.3}
\| q \|_{\ov{D}} \leq \beta^n \| p \|_{\ov{D}},
\end{equation}
with
\begin{equation} \label{1.4}
\beta := \exp \left( \frac{1}{\pi} \int_0^{2\pi/3} \log
\left(2 \cos \frac{t}{2}\right) dt \right).
\end{equation}
The constant $\beta=C_{\ov D}$ is asymptotically sharp, as $n \to \infty,$ and
it can also be expressed in a different way, using Mahler's measure. This
problem is of importance in designing algorithms for factoring
polynomials with integer coefficients over integers. We refer to \cite{Boy3}
and \cite{La} for more information on the connection with symbolic computations.

A further development related to (\ref{1.2}) for $E=[-a,a],\ a>0,$ was
suggested by P. B. Borwein in \cite{Bor} (see Theorems 2 and 5 there or see
Section 5.3 in \cite{BE}). In particular,  Borwein proved that if $\deg q = m$
then
\begin{equation} \label{1.5}
|q(-a)| \leq \| p \|_{[-a,a]} a^{m-n} 2^{n -1} \prod_{k =1}^{m} \left(1+\cos
\frac{2k -1}{2n} \pi \right),
\end{equation}
where the bound is attained for a monic Chebyshev polynomial of degree $n$ on
$[-a,a]$ and a factor $q$. He also showed that, for $E=[-2,2]$, the constant in
the above inequality satisfies
\begin{eqnarray*}
&&\limsup_{n \to \infty} \left(2^{m-1} \prod_{k =1}^{m} \left( 1 + \cos \frac{2k -1}{2n}
\pi \right)\right)^{1/n} \\ &\le& \lim_{n \to \infty} \left(2^{[2n/3]-1} \prod_{k =1}^{[2n/3]}
\left( 1 + \cos \frac{2k -1}{2n} \pi \right)\right)^{1/n} \\ &=& \exp \left( \int_0^{2/3}
\log \left(2 + 2\cos{\pi x}\right) dx \right) = 1.9081 \ldots,
\end{eqnarray*}
which hints that
\begin{equation} \label{1.6}
C_{[-2,2]} = \exp \left( \int_0^{2/3} \log \left(2 + 2\cos{\pi x}\right) dx
\right) = 1.9081 \ldots.
\end{equation}
We find the asymptotically best constant $C_E$ in (\ref{1.2}) for a rather
arbitrary compact set $E$. The general result is then applied to the cases of a
disk and a line segment, so that we recover (\ref{1.3})-(\ref{1.4}) and confirm
(\ref{1.6}).

\section{Results}

Our solution of the above problem is based on certain ideas from the
logarithmic potential theory (cf. \cite{Ra} or \cite{Ts}). Let ${\rm cap}(E)$ be the {\it
logarithmic capacity} of a compact set $E \subset {\C}$.   For $E$ with ${\rm
cap} (E) > 0$, denote the {\it equilibrium measure} of  $E$ (in the sense of
the logarithmic potential theory) by $\mu_E$. We remark that $\mu_E$ is a
positive unit Borel measure supported on $E$, ${\rm supp}\, \mu_E \subset E$
(see \cite[p. 55]{Ts}).

\begin{thm} \label{thm1}
Let $E \subset {\C}$ be a compact set, ${\rm cap} (E) >0$.  Then the best
constant $C_E$ in {\rm (\ref{1.2})} is given by
\begin{equation} \label{2.1}
C_E = \frac{\dis\max_{u \in \partial E} \exp\left( \int_{|z-u| \ge 1} \log
|z-u| d \mu_E (z)\right)}{{\rm cap} (E)}.
\end{equation}
Furthermore, if $E$ is regular then
\begin{equation} \label{2.2}
C_E = \max_{u \in \partial E} \exp\left(\dis -\int_{|z-u| \le 1} \log |z-u| d
\mu_E (z)\right).
\end{equation}
\end{thm}
The above notion of regularity is to be understood in the sense of the exterior
Dirichlet problem (cf. \cite[p. 7]{Ts}).  Note that the condition ${\rm cap}
(E) >0$ is usually satisfied for all applications, as it only fails for very
{\it thin} sets (see \cite[pp. 63-66]{Ts}), e.g., finite sets in the plane. But
if $E$ consists of finitely many points then the inequality (\ref{1.2}) cannot
be true for a polynomial $p(z)$ with zeros at every point of $E$ and for its
linear factors $q(z)$.  On the other hand, Theorem \ref{thm1} is applicable to
any compact set with a connected component consisting of more than one point
(cf. \cite[p. 56]{Ts}).

One can readily see from (\ref{1.2}) or (\ref{2.1}) that the best constant
$C_E$ is invariant under the rigid motions of the set $E$ in the plane.
Therefore we consider applications of Theorem \ref{thm1} to the family of disks
$D_r := \{z:|z|<r\}$, which are centered at the origin, and to the family of
segments $[-a,a],\ a>0.$

\begin{figure}[htb]
  \centerline{\psfig{figure=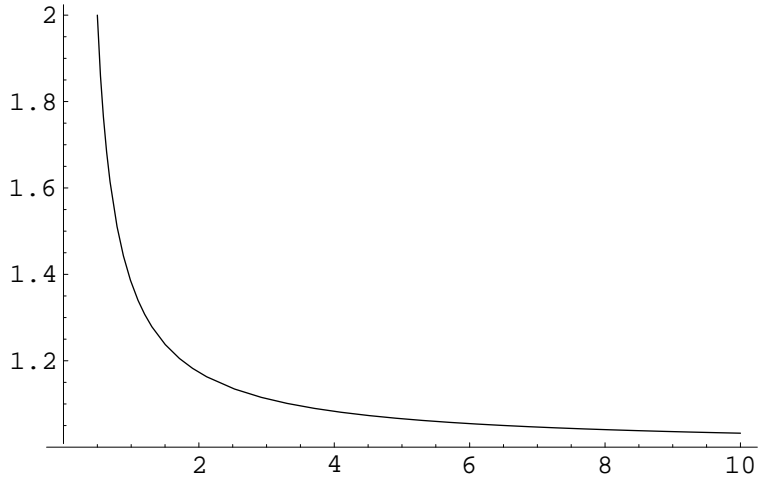,height=2.5in}}
  \caption{$C_{\ov D_r}$ as a function of $r$.}
  \label{fig1}
\end{figure}

\begin{cor} \label{cor1}
Let $D_r$ be a disk of radius $r$. Then the best constant $C_{\ov D_r}$, for
$E=\ov{D_r}$, is given by
\begin{equation} \label{2.3}
C_{\ov D_r} = \left \{ \begin{array}{l} \dis \frac{1}{r}, \quad 0<r \le 1/2, \\
\\ \dis \frac{1}{r} \exp \left(\frac{1}{\pi} \dis \int_0^{\pi-2\arcsin\frac{1}{2r}} \log
\left(2 r \cos\frac{x}{2}\right)\, dx \right), \quad  r>1/2.
\end{array} \right.
\end{equation}
\end{cor}
Note that (\ref{1.3})-(\ref{1.4}) immediately follow from (\ref{2.3}) for
$r=1.$ The graph of $C_{\ov D_r}$, as a function of $r$, is in Figure \ref{fig1}.

\begin{figure}[htb]
  \centerline{\psfig{figure=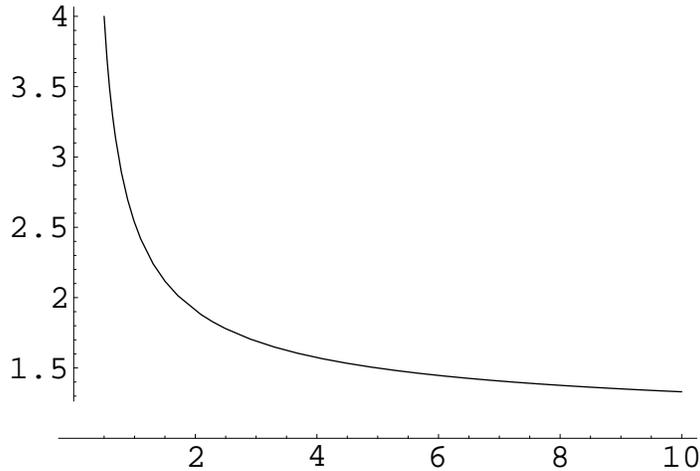,height=2.5in}}
  \caption{$C_{[-a,a]}$ as a function of $a$.}
  \label{fig2}
\end{figure}

\begin{cor} \label{cor2}
If $E=[-a,a],\ a>0$, then
\begin{equation} \label{2.4}
C_{[-a,a]} = \left \{ \begin{array}{l} \dis \frac{2}{a}, \quad 0< a \le 1/2,
\\ \\ \dis \frac{2}{a} \exp \left(\dis \int_{1-a}^a \frac{\log(t+a)}{\pi
\sqrt{a^2-t^2}}\, dt \right), \quad  a>1/2.
\end{array} \right.
\end{equation}
\end{cor}
Observe that (\ref{2.4}), with $a=2$, implies (\ref{1.6}) by the change of
variable $t=2\cos{\pi x}.$ We include the graph of $C_{[-a,a]}$, as a function
of $a$, in Figure \ref{fig2}.

We now state two general consequences of Theorem \ref{thm1}. They explain
some interesting features of $C_E$, which the reader may have noticed in Corollaries
\ref{cor1} and \ref{cor2}. Let
$$ {\rm diam}(E) := \max_{z,\zeta \in E} |z-\zeta|$$
be the Euclidean diameter of $E$.
\begin{cor} \label{cor3}
Suppose that ${\rm cap}(E)>0.$ If ${\rm diam}(E) \le 1$ then
\begin{equation} \label{2.5}
C_E = \frac{1}{{\rm cap}(E)}.
\end{equation}
\end{cor}
It is well known that cap$(D_r)=r$ and cap$([-a,a])=a/2$ (see \cite[p.
135]{Ra}), which clarifies the first lines of (\ref{2.3}) and (\ref{2.4}) by
(\ref{2.5}).

The next Corollary shows how the constant $C_E$ behaves under dilations of
the set $E$. Let $\alpha E$ be the dilation of $E$ with a factor $\alpha>0.$
\begin{cor} \label{cor4}
If $E$ is regular then
\begin{equation} \label{2.6}
\lim_{\alpha \to +\infty} C_{\alpha E} = 1.
\end{equation}
\end{cor}
Thus Figures \ref{fig1} and \ref{fig2} clearly illustrate (\ref{2.6}).

\medskip
We remark that one can deduce inequalities of the type (\ref{1.2}), for various
$L_p$ norms, from Theorem \ref{thm1}, by using relations between $L_p$ and
$L_{\infty}$ norms of polynomials on $E$ (see, e.g., \cite{Pr3}).

\section{Proofs}

\begin{proof}[Proof of Theorem \ref{thm1}]
The proof of this result is based on the ideas of \cite{Boy1} and \cite{Pr1}.
For $u \in \C$, consider a function
$$ \rho_u(z):=\max (|z-u|,1), \quad z \in \C.$$
One can immediately see that $\log \rho_u(z)$ is a subharmonic function in $z \in \C$,
which has the following integral representation (see \cite[p. 29]{Ra}):
\begin{equation} \label{3.1}
\log \rho_u(z) = \int \log |z-t|\, d\lambda_u(t), \quad z \in \C,
\end{equation}
where $d\lambda_u(u+e^{i\theta})=d\theta/(2\pi)$ is the normalized angular
measure on $|t-u|=1$.

Let $u \in \partial E$ be such that
$$ \|q\|_E=|q(u)|.$$
If $z_k,\ k =1, \ldots ,m,$ are the zeros of $q(z)$, counted according to multiplicities, then
\begin{eqnarray} \label{3.2}
\log \|q\|_E &=& \sum_{k =1}^m  \log |u-z_k| \leq \sum_{k =1}^m \log \rho_u(z_k)  \nonumber \\
&=& \sum_{k =1}^m \int \log |z_k -t|\, d\lambda_u(t) = \int \log |p(t)|\, d\lambda_u(t),
\end{eqnarray}
by (\ref{3.1}).

We use the well known Bernstein-Walsh lemma about the growth of a polynomial
outside of the set $E$ (see \cite[p. 156]{Ra}, for example):\\
Let $E \subset {\C}$ be a compact set, ${\rm cap}(E) >0$, with the unbounded component of
$\ov{{\C}} \setminus E$ denoted by $\Omega$.  Then, for any polynomial $p(z)$ of degree $n$,
we have
\begin{equation} \label{3.3}
|p (z)| \leq \| p \|_E \ e^{ng_{\Omega} (z, \infty )}, \quad z \in {\C},
\end{equation}
where $g_{\Omega} (z, \infty )$ is the Green function of $\Omega$, with pole at
$\infty$. The following representation for $g_{\Omega}(z,\infty)$ is found in Theorem III.37
of \cite[p. 82]{Ts}).
\begin{equation} \label{3.4}
g_{\Omega}(z,\infty)=\log \frac{1}{{\rm cap}(E)}+\int \log|z-t|\, d\mu_E(t), \quad z \in {\C}.
\end{equation}

It follows from (\ref{3.1})-(\ref{3.4}) and Fubini's theorem that
\begin{eqnarray*}
\frac{1}{n} \log \frac{\|q\|_E}{\| p\|_E} &\leq&
\int \log \frac{|p(t)|^{1/n}}{\| p\|_E^{1/n}}\, d \lambda_u(t)
\leq \int g_{\Omega}(t,\infty)\, d \lambda_u(t) \\ &=&
\log \frac{1}{{\rm cap}(E)} + \int \int \log|z-t|\, d\lambda_u(t) d\mu_E(z)\\
&=& \log \frac{1}{{\rm cap}(E)} + \int \log \rho_u(z)\, d\mu_E(z).
\end{eqnarray*}
Using the definition of $\rho_u(z)$, we obtain from the above estimate that
\begin{eqnarray*}
\| q \|_E &\leq& \left( \frac{\dis \max_{u \in \partial E}\exp\left(\int\log\rho_u(z)\, d\mu_E (z)\right)}
{{\rm cap}(E)} \right)^n \| p \|_E \\ &=& \left(\frac{\dis\max_{u \in \partial E} \exp\left(
 \int_{|z-u| \ge 1} \log |z-u|\, d \mu_E (z)\right)}{{\rm cap} (E)}\right)^n \| p \|_E.
\end{eqnarray*}
Hence
\begin{equation} \label{3.5}
C_E \le \frac{\dis\max_{u \in \partial E} \exp\left( \int_{|z-u| \ge 1} \log
|z-u|\, d \mu_E (z)\right)}{{\rm cap} (E)}.
\end{equation}

In order to prove the inequality opposite to (\ref{3.5}), we consider the
$n$-th Fekete points $\{ a_{k,n} \}_{k =1}^n$ for the set $E$ (cf. \cite[p. 152]{Ra}).
Let
$$ p_n (z) := \prod_{k =1}^n (z - a_{k,n})$$
be the Fekete polynomial of degree $n$. Define the normalized counting measures
on the Fekete points by
$$ \tau_n := \frac{1}{n} \sum_{k =1}^n \delta_{a_{k,n}}, \quad n \in {\N}.$$
It is known that (see Theorems 5.5.4 and 5.5.2 in \cite[pp. 153-155]{Ra})
\begin{equation} \label{3.6}
\lim_{n \rightarrow \infty} \| p_n \|_E^{1/n} = {\rm cap} (E).
\end{equation}
Furthermore, we have the following weak* convergence of counting measures
(cf. \cite[p. 159]{Ra}):
\begin{equation} \label{3.7}
\tau_n \stackrel{*}{\rightarrow} \mu_E, \quad \mbox{ as } n \rightarrow \infty.
\end{equation}
Let $u \in \partial E$ be a point, where the maximum on the right hand side of
(\ref{3.5}) is attained. Define the factor $q_n(z)$ for $p_n(z)$, with zeros being
the $n$-th Fekete points satisfying $|a_{k,n}-u| \ge 1$. Then we have by
(\ref{3.7}) that
\begin{eqnarray*}
\lim_{n \rightarrow \infty} \| q_n \|_E^{1/n} &\ge& \lim_{n \rightarrow \infty} |q_n(u)|^{1/n} =
\lim_{n \rightarrow \infty} \exp\left( \frac{1}{n} \sum_{|a_{k,n}-u| \ge 1}\log|u-a_{k,n}| \right)  \\
&=& \exp \left( \lim_{n \rightarrow \infty} \int_{|z-u| \ge 1} \log|u-z|\, d \tau_n (z)
\right) \\ &=& \exp \left( \int_{|z-u| \ge 1} \log|u-z| d \mu_E (z) \right).
\end{eqnarray*}
Combining the above inequality with (\ref{3.6}) and the definition of $C_E$, we
obtain that $$C_E \geq \lim_{n \rightarrow \infty} \frac{\|q_n\|_E^{1/n}}{\|
p_n \|_E^{1/n}} \ge \frac{\dis \exp\left(\int_{|z-u| \ge 1} \log |z-u|\, d
\mu_E (z)\right)}{{\rm cap} (E)}.$$ This shows that (\ref{2.1}) holds true.
Moreover, if $u \in \partial E$ is a regular point for $\Omega$, then we obtain
by Theorem III.36 of \cite[p. 82]{Ts}) and (\ref{3.4}) that $$ \log
\frac{1}{{\rm cap}(E)}+\int \log|u-t|\, d\mu_E(t) = g_{\Omega}(u,\infty) = 0.$$
Hence $$ \log \frac{1}{{\rm cap}(E)}+\int_{|z-u| \ge 1} \log|u-t|\, d\mu_E(t) =
-\int_{|z-u| \le 1} \log|u-t|\, d\mu_E(t), $$ which implies (\ref{2.2}) by
(\ref{2.1}).
\end{proof}

\begin{proof}[Proof of Corollary \ref{cor1}]
It is well known \cite[p. 84]{Ts} that cap$(\ov {D_r})=r$ and $d \mu_{\ov
{D_r}}(re^{i\theta}) = d\theta/(2\pi),$ where $d\theta$ is the angular measure on $\partial
{D_r}$. If $r \in (0,1/2]$ then the numerator of (\ref{2.1}) is equal to 1, so
that
$$ C_{\ov {D_r}} = \frac{1}{r}, \quad 0<r\le1/2. $$
Assume that $r>1/2.$ We set $z=re^{i\theta}$ and let $u_0=re^{i\theta_0}$ be a
point where the maximum in (\ref{2.1}) is attained. On writing
$$ |z-u_0|=2r\left|\sin\frac{\theta-\theta_0}{2}\right|,$$
we obtain that
\begin{eqnarray*}
C_{\ov {D_r}} &=& \frac{1}{r}\exp\left(\frac{1}{2\pi}\int_{\theta_0+2\arcsin\frac{1}{2r}}
^{2\pi+\theta_0-2\arcsin\frac{1}{2r}}
\log\left|2r\sin\frac{\theta-\theta_0}{2}\right|\, d\theta\right) \\ &=&
\frac{1}{r} \exp \left(\frac{1}{2\pi} \dis \int_{2\arcsin\frac{1}{2r}-\pi}^{\pi-2\arcsin\frac{1}{2r}}
\log \left(2 r \cos\frac{x}{2}\right) dx \right) \\ &=&
\frac{1}{r} \exp \left(\frac{1}{\pi} \dis \int_0^{\pi-2\arcsin\frac{1}{2r}} \log
\left(2 r \cos\frac{x}{2}\right) dx \right),
\end{eqnarray*}
by the change of variable $\theta-\theta_0=\pi-x.$
\end{proof}

\begin{proof}[Proof of Corollary \ref{cor2}]
Recall that cap$([-a,a])=a/2$ (see \cite[p. 84]{Ts}) and
$$ d\, \mu_{[-a,a]}(t)=\frac{dt}{\pi\sqrt{a^2-t^2}}, \quad t \in [-a,a].$$
It follows from (\ref{2.1}) that
\begin{equation} \label{3.8}
C_{[-a,a]} = \frac{2}{a}\, \exp \left(\dis\max_{u \in [-a,a]} \int_{[-a,a]\setminus(u-1,u+1)}
\frac{\log|t-u|} {\pi \sqrt{a^2-t^2}}\, dt \right).
\end{equation}
If $a \in (0,1/2]$ then the integral in (\ref{3.8}) obviously vanishes, so that
$C_{[-a,a]} = 2/a$. For $a>1/2$, let
\begin{equation} \label{3.9}
f(u):= \int_{[-a,a]\setminus(u-1,u+1)} \frac{\log|t-u|} {\pi \sqrt{a^2-t^2}}\, dt.
\end{equation}
One can easily see from (\ref{3.9}) that
$$ f'(u) = \int_{u+1}^a \frac{dt} {\pi(u-t) \sqrt{a^2-t^2}} < 0, \quad u \in
[-a,1-a],$$
and
$$ f'(u) = \int_{-a}^{u-1} \frac{dt} {\pi(u-t) \sqrt{a^2-t^2}} > 0, \quad u \in
[a-1,a].$$
However, if $u \in (1-a,a-1)$ then
$$f'(u) = \int_{u+1}^a \frac{dt} {\pi(u-t) \sqrt{a^2-t^2}} +
\int_{-a}^{u-1} \frac{dt} {\pi(u-t) \sqrt{a^2-t^2}}.$$
It is not difficult to verify directly that
$$\int \frac{dt} {\pi(u-t) \sqrt{a^2-t^2}} = \frac{1}{\pi\sqrt{a^2-u^2}} \log
\left|\frac{a^2-ut+\sqrt{a^2-t^2}\,\sqrt{a^2-u^2}}{t-u}\right| + C, $$
which implies that
$$ f'(u)= \frac{1}{\pi\sqrt{a^2-u^2}} \log \left( \frac{a^2-u^2+u+\sqrt{a^2-(u-1)^2}\,\sqrt{a^2-u^2}}
{a^2-u^2-u+\sqrt{a^2-(u+1)^2}\,\sqrt{a^2-u^2}} \right),  $$
for $u \in (1-a,a-1).$ Hence
$$ f'(u)<0,\ u \in (1-a,0), \quad \mbox{ and } \quad f'(u)>0,\ u \in (0,a-1). $$
Collecting all facts, we obtain that the maximum for $f(u)$ on $[-a,a]$ is
attained at the endpoints $u=a$ and $u=-a$, and it is equal to
$$\max_{u \in [-a,a]} f(u) = \int_{1-a}^a \frac{\log(t+a)}{\pi \sqrt{a^2-t^2}}\, dt.$$
Thus (\ref{2.3}) follows from (\ref{3.8}) and the above equation.
\end{proof}

\begin{proof}[Proof of Corollary \ref{cor3}]
Note that the numerator of (\ref{2.1}) is equal to 1, because $|z-u| \le 1,\
\forall\, z \in E,\ \forall\, u \in \partial E$. Thus (\ref{2.5}) follows
immediately.
\end{proof}

\begin{proof}[Proof of Corollary \ref{cor4}]
Observe that $C_E \ge 1$ for any $E \in \C$, so that $C_{\alpha E} \ge 1$.
Since $E$ is regular, we use the representation for $C_E$ in (\ref{2.2}). Let
$T:E \to \alpha E$ be the dilation mapping. Then $|Tz-Tu|=\alpha|z-u|,\ z,u \in E,$
and $d\mu_{\alpha E}(Tz)=d\mu_E(z)$. This gives that
\begin{eqnarray*}
C_{\alpha E} &=& \max_{Tu \in \partial (\alpha E)} \exp\left(\dis -\int_{|Tz-Tu| \le 1}
\log |Tz-Tu|\, d\mu_{\alpha E} (Tz)\right) \\ &=&
\max_{u \in \partial E} \exp\left(\dis -\int_{|z-u| \le 1/\alpha}
\log (\alpha |z-u|)\, d\mu_{E} (z)\right) \\ &=&
\max_{u \in \partial E} \exp\left(\dis -\mu_E(\ov{D_{1/\alpha}(u)})\log\alpha-
\int_{|z-u| \le 1/\alpha} \log |z-u|\, d\mu_{E} (z)\right) \\ &<&
\max_{u \in \partial E} \exp\left(\dis -\int_{|z-u| \le 1/\alpha} \log |z-u|\, d\mu_{E}
(z)\right),
\end{eqnarray*}
where $\alpha \ge 1$. Using the
absolute continuity of the integral, we have that
$$ \lim_{\alpha \to +\infty} \int_{|z-u| \le 1/\alpha} \log |z-u|\, d\mu_{E} (z)
=0, $$
which implies (\ref{2.6}).
\end{proof}


\end{document}